\documentclass{amsart}

\title[Weak coherence and finite decomposition complexity]
{Weak coherence and the \textit{K}-theory of groups with finite decomposition complexity}
\author[Boris Goldfarb]{Boris Goldfarb}
\address{Department of Mathematics and Statistics, SUNY, Albany, NY 12222}
\email{bgoldfarb@albany.edu}
\urladdr{http://www.albany.edu/~goldfarb/}
\date{\today}

\usepackage{amsmath,amsthm,amssymb,amsxtra,amscd,enumerate,mathrsfs,graphicx,pifont}
\usepackage[all,cmtip]{xy}
\SelectTips{cm}{10}\UseTips 

\theoremstyle{plain}

\theoremstyle{definition}           

{\leqno{\normalshape ({\thedummy})} $$}
\numberwithin{equation}{dummy}

\theoremstyle{plain}

\theoremstyle{plain}
\newtheorem{Thm}{Theorem}[section]

\newtheorem{Cor}[Thm]{Corollary}
\newtheorem{Lem}[Thm]{Lemma}
\newtheorem{Prop}[Thm]{Proposition}

\theoremstyle{definition}
\newtheorem{Def}[Thm]{Definition}

\newtheorem{Ex}[Thm]{Example}
\newtheorem{Rem}[Thm]{Remark}

\theoremstyle{remark}
\newtheorem{Not}[Thm]{Notation}

\DeclareMathOperator{\h}{\mathit{h}}
         \newcommand{\hlf}{\sideset{}{^{\textit{lf}}}\h}

\DeclareMathOperator{\im}{im}

\DeclareMathOperator{\K}{\mathit{K}}
         \newcommand{\Knc}{\K^{-\infty}}

\DeclareMathOperator{\Mod}{\mathbf{Mod}}

\DeclareMathOperator{\Wh}{\textit{Wh}}

\newcommand{\pull}
{\!\!\! -\!\!\! -\!\!\! -\!\!\!}
\DeclareMathOperator*{\hocolimprep}{hocolim}
\newcommand{\hocolim}[1]%
{\hocolimprep_{\substack{- \pull \rightarrow \\ #1}} \, }
\DeclareMathOperator*{\colimprep}{colim}
\newcommand{\colim}[1]%
{\colimprep_{\substack{- \pull \rightarrow \\ #1}} \, }

\newcommand{\define}[1]{\textit{#1}}

\providecommand{\bysame}{\makebox[3em]{\hrulefill}\thinspace}

\newif\ifShowLabels
\ShowLabelstrue
\newcommand{\TeXref}[1]{
\marginpar{\scriptsize \texttt{#1}}}

\newtheoremstyle{freestylethm}{6pt}{6pt}{\itshape}{}%
                {\bfseries}{}{.5em}{\thmnote{#3}}
\theoremstyle{freestylethm}
\newtheorem*{varthm}{}

\newcommand{\SecRef}[2]{\section{#1}\label{S:#2}%
\ifShowLabels \TeXref{{S:#2}} \fi}

\newcommand{\refT}[1]{\textup{\ref{T:#1}}}

\newcommand{\refC}[1]{\textup{\ref{C:#1}}}

\newcommand{\refP}[1]{\textup{\ref{P:#1}}}

\newenvironment{ThmRef}[1]%
{ \begin{Thm} \label{T:#1}
\ifShowLabels \TeXref{T:#1} \fi }%
{ \end{Thm} }
\newenvironment{DefRef}[1]%
{ \begin{Def} \label{D:#1}
\ifShowLabels \TeXref{D:#1} \fi }%
{ \end{Def} }
{ \begin{Lem} \label{L:#1}
\ifShowLabels \TeXref{L:#1} \fi }%
{ \end{Lem} }
\newenvironment{CorRef}[1]%
{ \begin{Cor} \label{C:#1}
\ifShowLabels \TeXref{C:#1} \fi }%
{ \end{Cor} }
{ \begin{Rem} \label{R:#1}
\ifShowLabels \TeXref{R:#1} \fi }%
{ \end{Rem} }
\newenvironment{PropRef}[1]%
{ \begin{Prop} \label{P:#1}
\ifShowLabels \TeXref{P:#1} \fi }%
{ \end{Prop} }
\newenvironment{ExRef}[1]%
{ \begin{Ex} \label{E:#1}
\ifShowLabels \TeXref{E:#1} \fi }%
{ \end{Ex} }
\newenvironment{NotRef}[1]%
{ \begin{Not} \label{N:#1}
\ifShowLabels \TeXref{N:#1} \fi }%
{ \end{Not} }

\newenvironment{ThmRefName}[2]%
{ \begin{Thm} [#2]\label{T:#1}
\ifShowLabels \TeXref{T:#1} \fi }%
{ \end{Thm} }
\newenvironment{DefRefName}[2]%
{ \begin{Def} [#2]\label{D:#1}
\ifShowLabels \TeXref{D:#1} \fi }%
{ \end{Def} }
{ \begin{Lem} [#2]\label{L:#1}
\ifShowLabels \TeXref{L:#1} \fi }%
{ \end{Lem} }
{ \begin{Cor} [#2]\label{C:#1}
\ifShowLabels \TeXref{C:#1} \fi }%
{ \end{Cor} }
{ \begin{Rem} [#2]\label{R:#1}
\ifShowLabels \TeXref{R:#1} \fi }%
{ \end{Rem} }
{ \begin{Prop} [#2]\label{P:#1}
\ifShowLabels \TeXref{P:#1} \fi }%
{ \end{Prop} }
{ \begin{Ex} [#2]\label{E:#1}
\ifShowLabels \TeXref{E:#1} \fi  }%
{ \end{Ex} }

\ShowLabelsfalse

\begin{document}

\begin{abstract}
The weak regular coherence is a coarse property of a finitely generated group $\Gamma$.
It was introduced by G.~Carlsson and this author to play the role of a weakening of Waldhausen's regular coherence as
part of computation of the integral $K$-theoretic assembly map.

A new class of metric spaces (sFDC) was introduced recently by A.~Dranishnikov and M.~Zarichnyi.
This class includes most notably the spaces with finite decomposition complexity (FDC) studied by E.~Guentner, D.~Ramras, R.~Tessera, and G.~Yu.
The main theorem of this paper shows that a group that has finite $K(\Gamma,1)$ and sFDC is weakly regular coherent.

As a consequence, the integral $K$-theoretic assembly maps are isomorphisms in all dimensions for
any group that has finite $K(\Gamma,1)$ and FDC.
In particular, the Whitehead group $\Wh (\Gamma)$ is trivial for such groups.
\end{abstract}

\maketitle

\section*{Introduction}

Given a topological manifold, there are many consequences of the fact that the fundamental group $\Gamma$ has finite asymptotic dimension (FAD).  It is known that the Novikov conjecture is true for such manifolds from the pioneering work of Yu \cite{gY:98}.  The Novikov conjecture in integral $K$-theory was proved for groups with FAD by Bartels \cite{aB:03}.  Further, the integral Borel conjecture for groups with finite $K(\pi,1)$ and FAD was verified by G.~Carlsson and the author \cite{gCbG:03,gCbG:04,gCbG:11,gCbG:13}.  The latter fact, in particular, has the consequence that the Whitehead group of such $\Gamma$ is trivial.

The ultimate goal of this paper is to extend these results further to groups that have finite decomposition complexity (FDC) as defined by Guentner--Tessera--Yu \cite{eGrTgY:12}.
The work of Carlsson--Goldfarb shows that the Borel isomorphism conjecture follows from the combination of (1) the coarse integral Borel conjecture and (2) the weak regular coherence property.
Both of these conditions are reviewed in the paper.
The first condition has already been verified for groups with a finite $K(\Gamma,1)$ and
finite decomposition complexity by Ramras--Tessera--Yu \cite{dRrTgY:11}.

The following are our main results.
The class $\mathrm{LH} \mathcal{F}$ of locally hierarchically decomposable groups was introduced by P.~Kropholler in \cite{pK:93}, the other terms are defined further in the paper.

\begin{varthm}[Main Theorem.]
Suppose $\Gamma$ is a torsion-free group from Kropholler's class $\mathrm{LH} \mathcal{F}$ and has straight finite decomposition complexity.
Then $\Gamma$ is weakly regular coherent.
\end{varthm}

\begin{varthm}[Corollary.]
Suppose, more specifically,
$\Gamma$ has a finite $K(\Gamma,1)$ and has
finite decomposition complexity.
Let $R$ be a noetherian ring of finite global dimension.
Then the integral $K$-theoretic assembly map $A(\Gamma,R)$ is a weak homotopy equivalence.
As a consequence, the Whitehead group $\textit{Wh} (\Gamma)$ is trivial for such groups.
\end{varthm}

The groups with FDC or sFDC are intermediate between the class FAD and the groups with Yu's property A \cite{aDmZ:12}.  Some well-known examples of groups with infinite asymptotic dimension, such as the wreath product of two copies of the integers $\mathbb{Z}$, have FDC.
Further examples of FDC groups established in \cite{eGrTgY:12} are
all finitely generated subgroups of $GL_n (k)$, where $k$ is a field.
For other interesting infinite dimensional groups such as Thompson's group, Grave's group \cite{bG:03}, Gromov's random groups the answer is unknown.  It is also unknown where groups such as $\textrm{Out} (F_n)$ belong in this hierarchy.

I am grateful to A.~Dranishnikov for sharing the early version of the paper \cite{aDmZ:12} and related communications.
I also thank the referees for corrections and thorough comments that significantly improved the exposition.

\SecRef{Geometric preliminaries}{LKOR}

The following is a summary of relationships between some coarse properties of metric spaces.

The original definition of asymptotic dimension by M.~Gromov is a coarse analogue of the covering dimension of topological spaces.

\begin{DefRef}{ASDIM}
The \textit{asymptotic dimension} of a metric space $X$ is defined as
the smallest number $n$ such that for any $d > 0$ there is a
uniformly bounded cover $\mathcal{U}$ of $X$
such that any metric ball of radius $d$ in
$X$ meets no more than $n + 1$ elements of the cover $\mathcal{U}$.
If such number exists for $X$, one says that $X$ has \textit{finite asymptotic dimension}.
\end{DefRef}

A number of authors generalized FAD in several directions in recent years.
One such generalization is asymptotic property C defined in Dranishnikov \cite{aD:00}.

Given a subset $S$ of a metric space $X$, we will use the notation $S[b]$ for the $b$-enlargement of $S$, that is the subset $\{ x \in X \mid d(x,s) \le b \ \mathrm{for \ some} \ s \in S \}$.
So, in particular, the metric ball centered at $x$ with radius $r$ is denoted by $x[r]$.
Also, given a number $R > 0$, a collection of disjoint subsets $S_{\alpha}$ of $X$ is called $R$-\textit{disjoint} if $S_{\alpha} [R]$ is disjoint from the union $\bigcup_{\beta \ne \alpha} S_{\beta}$, for all $\alpha$.

\begin{DefRefName}{APC}{APC}
A metric space $X$ has the \textit{asymptotic property C}
if for every sequence of positive numbers $R_1 \le R_2 \le \ldots$ there exists a natural number $n$ and
uniformly bounded $R_i$-disjoint families $W_i$, $1 \le i \le n$, such that the union of all $n$ families is a covering of $X$.
\end{DefRefName}

On the other hand, the following is one of the equivalent definitions of the finite decomposition complexity from Guentner--Tessera--Yu.

Let $\mathcal{X}$ and $\mathcal{Y}$ be two families of metric spaces.  Let $R > 0$.
The family $\mathcal{X}$ is called $R$-\textit{decomposable over} $\mathcal{Y}$ if for any space $X$ in $\mathcal{X}$ there are collections of subsets $\{ U_{1,\alpha} \}_{\alpha \in A}$, $\{ U_{2,\beta} \}_{\beta \in B}$ such that
\[
X = \bigcup_{\stackrel{\scriptstyle i=1,2}{\scriptstyle \gamma = A \cup B}} U_{i,\gamma},
\]
each $U_{i,\gamma}$ is a member of the family $\mathcal{Y}$, and each of the collections $\{ U_{1,\alpha} \}$ and $\{ U_{2,\beta} \}$ is $R$-disjoint.
A family of metric spaces is called \textit{bounded} if there is a uniform bound on the diameters of the spaces in the family.

One of the equivalent definitions of finite decomposition complexity of the metric space $X$ is in terms of a winning strategy for the following game between two players.
In this description we can assume that the families of metric spaces that appear in the decompositions are families of metric subspaces of $X$.
In round number 1 the first player selects a number $R_1 > 0$, the second player has to select a family of metric spaces $\mathcal{Y}_1$ and an $R_1$-decomposition of $\{ X \}$ over $\mathcal{Y}_1$.
In each succeeding round number $i$ the first player selects a number $R_i > 0$, the second player has to select a family of metric spaces $\mathcal{Y}_i$ and an $R_i$-decomposition of $\mathcal{Y}_{i-1}$ over $\mathcal{Y}_i$.
The second player wins the game if for some finite value of $i$ the family $\mathcal{Y}_i$ is bounded.

\begin{DefRefName}{FDC}{FDC}
A metric space $X$ has \textit{finite decomposition complexity} if the second player possesses a winning strategy in every game played over $X$.
\end{DefRefName}

The following property was defined by A.~Dranishnikov and M.~Zarichnyi in \cite{aDmZ:12}.

\begin{DefRefName}{DZP}{sFDC}
A metric space $X$ has \textit{straight finite decomposition complexity} if, for any sequence $R_1 \le R_2 \le \ldots$
of positive numbers, there exists a finite sequence of metric families
$\mathcal{V}_1$, $\mathcal{V}_2$, \ldots, $\mathcal{V}_n$ such that $\{ X \}$ is $R_1$-decomposable over $\mathcal{V}_1$, $\mathcal{V}_1$ is $R_2$-
decomposable over $\mathcal{V}_2$, etc., and the family $\mathcal{V}_n$ is bounded.
\end{DefRefName}

The relationship between the four geometric conditions on a metric space is established in \cite{aDmZ:12} and is expressed by the diagram
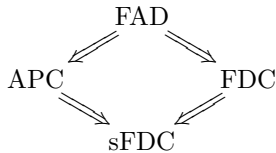
\[
\xymatrix@R.9pc@C.9pc{
&\textrm{FAD}
\ar@{=>}[ld]
\ar@{=>}[rd]
\\
\textrm{APC} \ar@{=>}[rd] &&\textrm{FDC} \ar@{=>}[ld]
\\
&\textrm{sFDC}
}
\]
All four classes satisfy Yu's property A.

In this paper we generalize \cite{gCbG:03} to a larger class of $R[\Gamma]$-modules by weakening the geometric condition on the group, promoting the main theorem from $\Gamma$ with FAD to $\Gamma$ with sFDC.

\SecRef{Filtered modules over a metric space}{Sheaves}

We start by recasting some of the controlled algebra from \cite{gCbG:11}.

Let $R$ be a left noetherian ring with unit and let $\Mod(R)$ be the category of left $R$-modules.
For a general set $X$, the power set $\mathcal{P}(X)$ partially ordered by inclusion is a category with subsets of $X$ as objects and unique morphisms $(S,T)$ when $S \subset T$.

\begin{DefRef}{MFIL}
An
\define{$X$-filtered $R$-module} is a
functor $F \colon \mathcal{P}(X) \rightarrow \Mod (R)$
where all structure maps $F(S,T)$ are monomorphisms and $F(\emptyset) = 0$.

We will abuse notation by referring to $F(X)$ as $F$.
For any covariant functor $F \colon \mathcal{P}(X) \to \Mod (R)$ with $F(\emptyset) = 0$ there exists an associated
$X$-filtered $R$-module $F_X$ given by
$F_X (S) = \im F (S,X)$.

Given a filtered module $F$ and an arbitrary submodule $F'$,
the submodule has the \textit{standard $X$-filtration} given by $F' (S) = F(S) \cap F'$.
For example, for each subset $T \subset X$, the submodule $F(T)$ of $F$ has the canonical filtration and gives a filtered
module $F_T$ with $F_T (S) = F(S) \cap F(T)$.
\end{DefRef}

We will assume that $X$ is a proper metric space in the sense that all closed metric balls in $X$ are compact.

\begin{DefRef}{KLIED}
Given a number $b \ge 0$,
an $R$-homomorphism $\varphi \colon F_1 \to F_2$ is called $b$-\textit{controlled} or simply \textit{controlled} if
the image
$\varphi (F_1 (S))$ is contained in the submodule $F_2 (S [b])$ for
all subsets $S \subset X$.

We say $\varphi$ is $b$-\textit{bicontrolled} or \textit{bicontrolled} if in addition to containments
\[
\varphi (F_1 (S)) \subset F_2 (S [b])
\]
as above, there are containments
\[
\varphi (F_1) \cap F_2 (S) \subset \varphi F_1 (S[b])
\]
for all subsets $S$ of $X$.

The $X$-filtered objects $F$ may be subject to the following constraints.
\begin{itemize}
\item $F$ is \textit{locally finite} if $F (V)$ is a finitely generated submodule of $F$ whenever $V$ is a bounded subset of $X$.
\item $F$ is $D$-\textit{lean} or simply \textit{lean} if there is a
number $D \ge 0$ such that we have
\[
F \left( U \right) \subset \sum_{x \in U} F \left( x [D]) \right)
\]
for any subset $U$ of $X$.
\item $F$ is $\delta$-\textit{split} or \textit{split} if there is a
number $\delta \ge 0$ such that we have
\[
F \left( U_1 \cup U_2 \right) \subset F(U_1 [\delta]) + F(U_2 [\delta])
\]
for any pair of subsets $U_1$, $U_2$ of $X$.
\item $F$ is $d$-\textit{insular} or \textit{insular} of there is a number
$d \ge 0$ such that
\[
F(U_1) \cap F(U_2) \subset F \big( U_1 [d] \cap U_2 [d] \big)
\]
for any pair of subsets $U_1$, $U_2$ of $X$.
\end{itemize}
\end{DefRef}

It is clear that a $D$-lean module is $D$-split.

If $F$ is $X$-filtered, the kernel $K$ of an arbitrary
homomorphism $f \colon F \to G$ inherits the
standard filtration
$K(S) = K \cap F(S)$.

\begin{PropRef}{TKOBBIS}
Given a bicontrolled epimorphism $f \colon F \to G$, suppose $F$ is split and $G$ is insular.
Then the kernel $K$ with the standard filtration is split.
\end{PropRef}

\begin{proof}
Suppose $F$ is $\delta$-split and $G$ is $d$-insular.  Let $i \colon K \to F$ be the inclusion of the kernel.
Given $z \in K (T \cup U)$, we have $i(z) \in F (T \cup U) \subset F (T[\delta]) + F (U[\delta])$,
so we can write accordingly $i(z) = y_1 + y_2$ where
$f(y_1) = - f(y_2)$.
Suppose the epimorphism $f$ is $b$-bicontrolled.
Since $G$ is $d$-insular,
\[
f(y_1) = - f(y_2) \in G \big( T[\delta + b] \big) \cap G \big( U[\delta + b] \big) \subset G \big( T[\delta + b+ d] \cap U[\delta + b + d] \big).
\]
Therefore we are able to find
\[
y \in F \big( T[\delta + 2b + d] \cap U[\delta + 2b + d] \big)
\]
such that $f(y) = f(y_1) = - f(y_2)$, because generally $(S \cap P)[b] \subset S[b] \cap P[b]$.
Thus
\[
i(z) = y_1 + y_2 = (y_1 - y) + (y_2 + y)
\]
and
\[
y_1 - y \in F \big( T[\delta + 2b + d] \big),
\quad
y_2 + y \in F \big( U[\delta + 2b + d] \big).
\]
Let $z_1 = i^{-1} (y_1 - y)$ and $z_2 = i^{-1} (y_2 + y)$, and we have $z = z_1 + z_2$ such that
\[
z_1 \in K \big( T[\delta + 2b + d] \big),
\quad
z_2 \in K \big( U[\delta + 2b + d] \big),
\]
so $K$ is $(\delta + 2b + d)$-split.
\end{proof}

\begin{PropRef}{CVBBTD}
Given a $b$-controlled homomorphism $f \colon F \to G$, suppose $F$ is $D$-lean and $G$ is $d$-insular.
Let $K$ be the kernel of $f$ and let $U$ be the union of a $(2D + 2b +2d)$-disjoint collection $\{ U_{\alpha} \}$ of subsets of $X$.
Then $K(U) \subset \sum K(U_{\alpha} [D])$.
\end{PropRef}

\begin{proof}
Since $F$ is $D$-lean, an arbitrary element $k \in K(U)$ can be written as a sum $k = \sum k_{\alpha}$, where $k_{\alpha} \in F(U_{\alpha} [D])$ and $k_{\alpha} =0$ for all but finitely many $\alpha$.  To check that $k_{\alpha} \in K$, consider $k'_{\alpha}= \sum_{\beta \ne \alpha} k_{\beta}$ and $U'_{\alpha} = \bigcup_{\beta \ne \alpha} U_{\beta}$.  Then $f(k_{\alpha}) \in G(U_{\alpha} [D + b])$ and $f(k'_{\alpha}) \in G(U'_{\alpha} [D + b])$.  So from the $d$-insularity assumption, $f(k_{\alpha}) = - f(k'_{\alpha}) \in G \big( U_{\alpha} [D + b + d] \cap U'_{\alpha} [D + b + d] \big) = G(\emptyset) = 0$.
\end{proof}

\begin{ThmRef}{CVBBCC}
Suppose $X$ is a metric space with sFDC.  Given a bicontrolled epimorphism $f \colon F \to G$, suppose $F$ is lean and $G$ is insular.
Then the kernel $K$ with the standard filtration is lean.
\end{ThmRef}

\begin{proof}
Let us assume that $f$ is $b$-bicontrolled, $F$ is $D$-lean, and $G$ is $d$-insular.
Since $F$ is $D$-split if it is $D$-lean, Proposition \refP{TKOBBIS} can be applied showing that $K$ is $D'$-split for $D' = D + 2b + d$.

We start with the sequence of numbers $R_1 \le R_2 \le \ldots \le R_k \le \ldots$ where we choose $R_k = 2kD' + 2kD +2b +2d$.
By the assumption, there is a sequence of metric families $\mathcal{Y}_0 = \{ X \}$, 
$\mathcal{Y}_1$, $\mathcal{Y}_2$, \ldots with the corresponding $R_i$-decompositions 
of $\mathcal{Y}_{i-1}$ over $\mathcal{Y}_i$ 
so that for some $n$ the diameter of every member of the family $\mathcal{Y}_n$ is bounded by the same number $M$.

To start the inductive argument, assume that $X = U_1 \cup U_2$ so that $U_1$ is the union of an $R_1$-disjoint family $\{ U_{1,\alpha} \}$, for $R_1 = 2D' + 2D+2b+2d$, and $U_2$ is the union of an $R_1$-disjoint family $\{ U_{2,\beta} \}$.
Then by Proposition \refP{TKOBBIS} each element $k$ in $K$ is the sum $k_{1} + k_{2}$ for some $k_{1} \in K(U_1 [D'])$ and $k_{2} \in K(U_2 [D'])$.
Since both $\{ U_{1,\alpha} [D'] \}$ and $\{ U_{2,\beta} [D'] \}$ are $(2D + 2b + 2d)$-disjoint, we can use Proposition \refP{CVBBTD} to write $k_{1} = \sum k_{1,\alpha}$, $k_{2} = \sum k_{2,\beta}$, where $k_{1,\alpha} \in K(U_{1,\alpha} [D'+D])$, $k_{2,\beta} \in K(U_{2,\beta} [D'+D])$ and $k_{1,\alpha} = k_{2,\beta} =0$ for all but finitely many $\alpha$ and $\beta$.
By the assumption, both $\{ U_{1,\alpha} \}$ and $\{ V_{1,\beta} \}$ are $R_2$-decomposable over $\mathcal{Y}_2$ for $R_2 = 4D' + 4D +2b + 2d$.  Therefore all $\{ U_{1,\alpha} [D'+D] \}$ and $\{ U_{2,\beta} [D'+D] \}$ are $(2D' + 2D +2b + 2d)$-decomposable over $\mathcal{Y}_2$.  
The decompositions can be continued inductively for the total of $n$ steps giving $k$ as the sum of kernel elements in $F(W[nD' + nD])$, where $W \in \mathcal{Y}_n$.  So $k$ is generated by the kernel elements with uniformly bounded support: if the diameters of sets in $\mathcal{Y}_n$ are bounded by $M$ then $K$ is $(M + nD' + nD)$-lean.
\end{proof}

We will view a finitely generated group $\Gamma$ as a metric space with the word metric generated from some choice of a finite generating set.  
The word-length metric on
$\Gamma$ is the path metric
induced from the condition that $d (\gamma, \gamma \omega) =
1$ whenever $\gamma \in \Gamma$ and $\omega \in \Omega$.
Then $\Gamma$ acts by left translations on itself, and each left translation is an isometry.

\begin{DefRef}{HUT}
Consider an $R[\Gamma]$-module $F$ and assume that $F$ is given
a locally finite $\Gamma$-filtration by $R$-submodules.
We will say that $F$ is $\Gamma$-\textit{equivariant}, or simply a $\Gamma$-\textit{module}, if $F(\gamma S) = \gamma F(S)$ for all choices of $\gamma \in \Gamma$ and $S \subset \Gamma$.
\end{DefRef}

The kernel, the cokernel, and the image of an $R[\Gamma]$-homomorphism $f \colon F \to G$ of $\Gamma$-modules are again $\Gamma$-modules.  If $K$ is the kernel of $f$, then the $\Gamma$-filtration of $K$ is given by $K(S) = K \cap F(S)$.
If $C$ is the cokernel, then $C(S)$ is the image of $G(S)/\left( f(F) \cap G(S) \right)$ under the canonical map to $C=G/f(F)$.

\begin{CorRef}{JKUF}
Suppose $R$ is a noetherian ring and $\Gamma$ is a finitely generated group with sFDC.  Given a bicontrolled $R[\Gamma]$-linear epimorphism $f \colon F \to G$ of $\Gamma$-modules, suppose $F$ is lean and $G$ is insular.
Then the kernel $K$ is finitely generated.
\end{CorRef}

\begin{proof}
The kernel is lean by Theorem \refT{CVBBCC}.  If $D$ is a leanness constant then $K$ is generated by $K \cap F(e[D])$, where $e$ is the identity element, as an $R[\Gamma]$-module.  Since $R$ is noetherian, the $R$-module $K \cap F(e[D])$ is finitely generated, so $K$ is finitely generated as an $R[\Gamma]$-module.
\end{proof}

\SecRef{Applications}{Appl}

\subsection{Weak regular coherence of sFDC groups}

Let $A$ be a ring with a unit.

\begin{DefRef}{HUYIV}
A left $A$-module is said to be of type $\mathrm{FP}_{\infty}$ if
it has a resolution by finitely generated
projective $A$-modules.
The ring $A$ is called \textit{coherent} if every finitely presented $A$-module is of type $\mathrm{FP}_{\infty}$.
It is called \textit{regular coherent} if each resolution
can be chosen to be finite.

An $A$-module has \textit{finite projective dimension} if it has a finite resolution by finitely generated projective $A$-modules.
The ring $A$ has \textit{finite (global) dimension}
if every finitely generated $A$-module has a resolution of some fixed length by finitely generated projective $A$-modules.
\end{DefRef}

Waldhausen \cite{fW:78} discovered a
remarkable collection of discrete groups $\Gamma$ such that all
finitely presented modules over the group ring $\mathbb{Z}[\Gamma]$ are
regular coherent. It includes free groups, free abelian groups,
torsion-free one relator groups, fundamental groups of
submanifolds of the three-dimensional sphere, and their various
amalgamated products and HNN extensions and so, in particular, the
fundamental groups of submanifolds of the three-dimensional
sphere.  Waldhausen called this property of the group \textit{regular
coherence} and used it to compute the algebraic $K$-theory of
these groups. He asked if a weaker property of the group
ring would suffice in his argument,
see for example the paragraph after the proof of Theorem 11.2 in~\cite{fW:78}.
The regular coherence property for groups seems to be very special---simply constructing
individual finite dimensional modules over group rings is hard.
The \textit{weak regular coherence} property
defined in Carlsson--Goldfarb \cite{gCbG:03} turned out to play that role
in a different argument \cite{gCbG:13}
with the same purpose of computing the integral assembly map in $K$-theory.

Given a left $R[\Gamma]$-module $F$ with a finite generating set
$\Sigma$, it is also an $R$-module with the generating set
$\widetilde{\Sigma} = \{ \gamma \sigma \in F \mid \gamma \in \Gamma, \sigma \in \Sigma \}$.
Now one can associate to every subset $S$ of $\Gamma$
the left $R$-submodule $F(S)$ generated by $\gamma \sigma \in \widetilde{\Sigma}$ such that
$\gamma \in S$ and $\sigma \in \Sigma$.
This gives a functor $F \colon \mathcal{P}(\Gamma) \to \mathrm{Mod}_R (F)$, from
the power set of $\Gamma$ to the $R$-submodules of $F$ such that $F(\Gamma) = F$, $F(\emptyset) = 0$, and for a bounded subset
$T \subset \Gamma$, $F(T)$ is a finitely generated $R$-module.
This shows that $F$ with the given filtration, which will be denoted by $F_{\Sigma}$, is a $\Gamma$-filtered $R$-module.
It is easy to see that $F_{\Sigma}$ is $\Gamma$-equivariant, so $F_{\Sigma}$ is a $\Gamma$-module.
Clearly, $F_{\Sigma}$ is $0$-lean by design.

For the future reference, let us remark that a finitely generated free $R[\Gamma]$-module with a finite set of free generators $\Sigma$ can be given the structure of a $\Gamma$-module as above.  In this case $\widetilde{\Sigma} = \Gamma \times \Sigma$.
It is easy to see that in this case the $\Gamma$-filtration is $0$-lean and $0$-insular.

\begin{NotRef}{HVC}
In the rest of this section, $R$ will be a noetherian ring with unit.
\end{NotRef}

The following are elementary facts about $\Gamma$-modules.

\begin{PropRef}{LeanBC}
Given $R[\Gamma]$-modules $F$ and $G$, let $\Sigma$ and $\Delta$ be finite generating sets for respectively $F$ and $G$.

(1) Every $R[\Gamma]$-homomorphism $\phi \colon F_{\Sigma} \to G_{\Delta}$ between finitely generated $R[\Gamma]$-modules is controlled.

(2) Every controlled $R[\Gamma]$-epimorphism $\psi \colon F \to G$ between lean $\Gamma$-modules is bicontrolled.  Therefore, every
$R[\Gamma]$-epimorphism $\phi \colon F_{\Sigma} \to G_{\Delta}$ is bicontrolled.

(3) Whenever two choices of finite generating sets $\Sigma$ and $\Delta$ give two filtrations $F_{\Sigma}$ and $F_{\Delta}$ of the same $R[\Gamma]$-module $F$, the identity map $F_{\Sigma} \to F_{\Delta}$ is bicontrolled.

(4) The kernel of a controlled $R$-homomorphism $\psi \colon F \to G$, where $F$ is insular, is insular.

(5) The cokernel and the image of a bicontrolled $R$-homomorphism $\psi \colon F \to G$ of lean insular
$\Gamma$-modules is lean and insular.
\end{PropRef}

\begin{proof}
(1) Consider $x \in F_{\Sigma}(S)$, which is the $R$-submodule generated by $S \Sigma$, then
\[
x = \sum_{s, \sigma} r_{s, \sigma} s \sigma
\]
for a finite collection of pairs $s \in S$, $\sigma \in \Sigma$.
For the identity element $e$ in $\Gamma$, $F_{\Sigma} (\{ e \})$ is the $R$-submodule generated by $\Sigma$.
Depending on the choice of $\Delta$ there is a number $b \ge 0$ such that
$\phi F_{\Sigma} (\{ e \}) \subset G_{\Delta} (e[b])$.
Therefore,
\[
\phi (x) = \sum_{s, \sigma} r_{s, \sigma} \phi (s \sigma) =
\sum_{s, \sigma} r_{s, \sigma} s \phi (\sigma) \subset \sum_{s \in
S} s \, G_{\Delta} (e[b]) \subset G_{\Delta} (S[b])
\]
because the left translation action by any element $s \in S$ on $e[b]$
in $\Gamma$ is an isometry onto $s[b]$.

(2) Suppose $F$ is $D_F$-lean, $G$ is $D_G$-lean, and $\psi$ is $b$-controlled.
We claim that there is a number $c \ge 0$ such that 
$G(x[D_G]) \subset \psi F(x[c])$ for a given point $x \in \Gamma$.
$G (S) \subset \sum_{x \in S} G(x[D])$.
If so then the same is true for all $x$ in $\Gamma$ by equivariance.
Assuming the claim,
\[
G (S) \subset \sum_{x \in S} G(x[D_G]) \subset \sum_{x \in S} \psi F(x[c]) \subset \psi F(S[c]).
\]
This shows that $\psi$ is $B$-bicontrolled for $B = \max \{ b, c \}$.

To prove the claim, choose finitely many $R$-linear generators for $G (x[D_F])$, say $g_1, \ldots, g_t$.
Choose arbitrary elements $f_1, \ldots, f_t$ in $F$ such that $\psi (f_i) = g_i$.
Since $F$ is $D_F$-lean, each $f_i = \sum f_{ij}$ for some elements $f_{ij} \in F (x_j [D_F])$ for a finite collection of $x_j$ in $\Gamma$.
Pick a number $C \ge 0$ such that all $x_j \in x[C]$.
Now all $f_i \in F (x[C + D_F])$.
So $G(x[D_G]) \subset \psi F(x[C + D_F])$, and $c$ can be chosen to be $C + D_F$, as claimed.

(3) is a consequence of (1) and (2).

(4) Recall that the filtration of the kernel $K$ is given by $K(S) = K \cap F(S)$.  If $F$ is $d$-insular then $K$ is $d$-insular.

(5) Since the quotient homomorphism is $0$-bicontrolled, the cokernel $C$ is $D$-lean if $G$ is $D$-lean.  The insularity property of $C$ follows from insularity of $G$ and leanness of $F$.  For details we refer to Proposition 2.18 of \cite{gCbG:11}.
For the image $I$ of $\psi$, the insularity property follows from part (4).  If $\psi$ is $b$-bicontrolled and $F$ is $D$-lean then it is easy to see that $I$ is $(2b + D)$-lean.
\end{proof}

One has the following family of important examples.

\begin{ExRef}{ExInt}
A controlled idempotent homomorphism of a
filtered module is always bicontrolled.  Indeed, if
$\phi \colon F \to F$ is an idempotent so that $\phi^2 = \phi$
then the restriction of $\phi$ to its image is the identity.
Therefore $\phi f(X) \cap f(S) \subset \phi
f(S)$.
As a consequence of part (5) of Proposition \refP{LeanBC}, the cokernels and images of idempotents of finitely generated free $R[\Gamma]$-modules are
examples of lean insular $\Gamma$-modules.

There are many systematic constructions of non-free stably free $\mathbb{Z}[\Gamma]$-modules for infinite discrete groups $\Gamma$, cf.~section 17.4 of Johnson \cite{fJ:12}. 
Specifically, there are constructions for geometrically interesting groups by Dunwoody \cite{mD:72}, Berridge--Dunwoody \cite{pBmD:79} for the trefoil group $\langle x,y \mid x^3 = y^2 \rangle$ and similar groups in the work concerned with relation modules. The recent constructions of Harlander--Jensen \cite{jHjJ:06} apply to the Baumslag--Solitar group $G = \langle x,y \mid x y^2 x^{-1} = y^3 \rangle$.
It is also true that $G$ is generated by $x$ and $z=y^4$.
Let $K = R(G,x,z)$ be the kernel of the canonical map $F_2 \to G$ sending the two free generators to $x$ and $z$.
The abelianization $M(G,x,z)=K/[K,K]$ is a $\mathbb{Z}[G]$-module, where the action of $G$ is by conjugation, called the relation module associated with the generating set $\{ x,z \}$. Now Harlander and Jensen show that $M$ is non free but $M \oplus \mathbb{Z}[G] \cong \mathbb{Z}[G]^2$.
The projection onto $M$ is an interesting idempotent $e$ of $\mathbb{Z}[G]^2$.
Let $e_n \colon \mathbb{Z}_n [G]^2 \to \mathbb{Z}_n [G]^2$ be the reduction modulo $n > 1$. The image $M_n$ of $e_n$ is not flat over $\mathbb{Z} [G]$, so it is not projective.
However, $e_n$ is a bicontrolled endomorphism of a lean insular $\mathbb{Z}[G]$-module $\mathbb{Z}_n [G]^2$, so $M_n$ is lean and insular.  
\end{ExRef}

\begin{DefRef}{FPMADM}
An $R[\Gamma]$-module is \textit{finitely presented} if it is the cokernel
of a homomorphism, called a \textit{presentation}, between
finitely generated free $R[\Gamma]$-modules. If the homomorphism is
bicontrolled, we call the presentation
\textit{admissible}.
\end {DefRef}

\begin{PropRef}{LeanBC2}
Let $\Gamma$ be a finitely generated group with sFDC.
Every lean insular
$\Gamma$-module
has an admissible presentation.
\end{PropRef}

\begin{proof}
Suppose $F$ is insular and $D$-lean.
Since the $R$-submodule $F(e[D])$ is finitely generated, there is finitely generated free $R$-module $F_e$ and an $R$-linear epimorphism $\phi_e \colon F_e \to F(e[D])$.  One similarly has epimorphisms $\phi_{\gamma} \colon F_{\gamma} \to F(\gamma[D])$ using isomorphic copies $F_{\gamma}$ of $F_e$.
We define a new locally finite $\Gamma$-filtered module $F_0$ by assigning $F_0 (S) = \bigoplus_{\gamma \in S} F_{\gamma}$.  Then clearly $F_0$ is a lean insular $\Gamma$-module.
There is an $R[\Gamma]$-homomorphism $\phi_0 \colon F_0 \to F$ sending $F_0 (\gamma)$ onto $F(\gamma[D])$ via $\phi_{\gamma}$.
Viewed as an $R$-homomorphism, this $\phi_0$ is $D$-bicontrolled, so from Theorem \refT{CVBBCC} the kernel $K$ of $\phi_0$ is a lean $\Gamma$-module.
It is insular by part (4) of Proposition \refP{LeanBC}.
One similarly constructs a filtered module $F_1$ and a bicontrolled $R[\Gamma]$-homomorphism $\phi_1 \colon F_1 \to K$.
Since $K$ is finitely generated by Corollary \refC{JKUF}, $F_1$ can be chosen to be finitely generated as well.
The composition of $\phi_1$ and the inclusion of $K$ in $F_0$ gives a bicontrolled $\psi_1 \colon F_1 \to F_0$ with the cokernel $F$, as required.
\end{proof}

\begin{CorRef}{LeanBCcoh}
Let $\Gamma$ be a finitely generated group with sFDC.
Every lean insular
$\Gamma$-module is of type $\mathrm{FP}_{\infty}$.
\end{CorRef}

\begin{proof}
Applying the construction from the proof of Proposition \refP{LeanBC2},
one obtains a projective $R[\Gamma]$-resolution by finitely generated $R[\Gamma]$-modules
\[
\ldots \longrightarrow F_{n} \longrightarrow F_{n-1}
\longrightarrow \ldots \longrightarrow F_1 \longrightarrow F_0 \longrightarrow F
\longrightarrow 0
\]
In this resolution, each homomorphism $\psi_n \colon F_n \to F_{n-1}$
factors through an epimorphism $\phi_n$ onto a lean insular $\Gamma$-submodule $K_{n-1}$ of $F_{n-1}$.
\end{proof}

The following are appropriate weakenings of the coherence conditions for finitely generated groups.

\begin{DefRef}{WRC}
A finitely generated group $\Gamma$ is \textit{weakly coherent} if for an arbitrary noetherian ring $R$ with unit every $R[\Gamma]$-module with an admissible presentation is of type $\mathrm{FP}_{\infty}$.

We say $\Gamma$ is \textit{weakly regular coherent} if,
whenever $R$ is coherent and has finite global dimension, every
$R[\Gamma]$-module with an admissible presentation has finite
projective dimension.
\end{DefRef}

Recall that P.~Kropholler \cite{pK:93} defined the class of groups $\mathrm{LH} \mathcal{F}$, which includes in particular all groups with finite $K(\Gamma,1)$.

The following is the main result of the paper.

\begin{ThmRefName}{Main1}{Main Theorem}
Suppose $\Gamma$ is a finitely generated group with the sFDC property.
Then $\Gamma$ is weakly coherent.
If, in addition, $\Gamma$ is torsion-free and belongs to the class $\mathrm{LH} \mathcal{F}$ then $\Gamma$ is weakly regular coherent.
\end{ThmRefName}

\begin{proof}
The first statement follows directly from Corollary \refC{LeanBCcoh}.
For the second statement we
examine the resulting projective resolution of $F$ with an admissible presentation.
By Theorem A of Kropholler \cite{pK:99},
there is $d \ge 0$ such that the $d$-th kernel $K_d$ is isomorphic to a direct summand of a polyelementary module.
When $\Gamma$ is torsion-free, every elementary module is of the form $U \otimes_R R[\Gamma]$, where $U$ is projective over $R$.  It follows that all elementary modules are projective over $R[\Gamma]$.
So the polyelementary modules are also projective.
This makes $K_d$ a projective $R[\Gamma]$-module.
\end{proof}

It is worth noting that the integer $d$ in the proof may be strictly greater than the projective dimension of the ring $R$.
There is still a possibility that for the lean insular $\Gamma$-modules, $R$-flatness implies $R[\Gamma]$-flatness, and therefore $R[\Gamma]$-projectivity since all lean insular $\Gamma$-modules are finitely presented.  This is certainly not true in general due to phenomena related to Moore's conjecture, cf.~\cite{eAeM:11}.

\begin{CorRef}{LPDSB}
Let $R=\mathbb{Z}$, and let $F$ be a lean insular $\Gamma$-module.  Suppose $\Gamma$ has sFDC, and there is a finite $K(\Gamma,1)$ complex.
Then $F$ has finite projective dimension over $\mathbb{Z}[\Gamma]$.
\end{CorRef}

\subsection{The Whitehead group of FDC groups}

The motivating goal of this paper is to extend the techniques used in \cite{gCbG:03,gCbG:13} to larger geometric contexts.
The results of this paper in conjunction with the proof of the Novikov conjecture for the FDC groups in \cite{dRrTgY:11} prove the vanishing of the Whitehead groups for FDC groups with finite classifying spaces.

The lean insular $\Gamma$-modules over a noetherian ring $R$ form a Quillen exact category which contains the finitely generated projective $R[\Gamma]$-modules as an exact subcategory.  For details we refer to \cite{gCbG:11}.
It follows from the Resolution Theorem of Quillen and the Main Theorem \refT{Main1} that the inclusion induces the Cartan isomorphisms on the level of algebraic $K$-groups when $\Gamma$ has sFDC and a finite $K(\Gamma,1)$ and when $R$ is coherent and has finite global dimension.
This applies in particular in the geometrically important case of $R = \mathbb{Z}$.

It is known that the Whitehead group $\textit{Wh} (\Gamma)$ contains obstructions to some crucial constructions in the classification theory of manifolds with the fundamental group $\Gamma$.  Vanishing of the entire Whitehead group is a very much desired fact.  This group is the cokernel of Loday's assembly map $A_1 \colon H_1 (\Gamma, K(\mathbb{Z})) \to K_1 (\mathbb{Z} [\Gamma])$, see \cite{jlL:76}, p.~311.  Thus $\textit{Wh} (\Gamma) = 0$ if $A_1$ is an epimorphism.  A conjecture in algebraic $K$-theory predicts that $A_1$ is an isomorphism for all torsion-free groups $\Gamma$. It is answered positively in \cite{gCbG:13} for a class of groups satisfying some specific conditions.
In order to state the result in its proper generality, 
we note that the conjecture about $A_1$ is part of a general conjecture about the assembly map 
\[
A \colon H (B\Gamma, \Knc (R)) \longrightarrow \Knc (R[\Gamma])
\]
between the nonconnective spectra representing the homology of $\Gamma$ with coefficients in the $K$-theory of the ring $R$ and the $K$-theory of the group ring $R[\Gamma]$.
The general statement that has no known counterexample is that $A$ is a weak equivalence for a regular noetherian ring $R$ and any torsion-free group $\Gamma$. 
Now the main theorem of \cite{gCbG:13} shows that $A$ is an equivalence 
whenever
(1) the ring $R$ has finite global dimension,
(2) there is a model for the Eilenberg--MacLane space $K(\Gamma,1)$ which is a finite complex, 
(3) the $\Gamma$-equivariant coarse assembly map
\[
a \colon \hlf (\Gamma, \Knc (R)) \longrightarrow \Knc (\Gamma, R)
\]
of nonconnective spectra that model the locally finite homology of $\Gamma$ and the bounded $K$-theory of $\Gamma$ is an equivalence, and
(4) the group $\Gamma$ is weakly regular coherent, giving therefore the Cartan equivalence from the preceding paragraph. 
We refer the reader to the statements of the Main Theorem (the expanded version) on page 3 and the statement in section 5.7 in \cite{gCbG:13}.

It is known that (3) holds for groups satisfying (2) and having FDC by a slightly stronger result of Ramras--Tessera--Yu \cite{dRrTgY:11}.
Condition (4) holds for groups with sFDC and a finite $K(\Gamma,1)$ by the main result of this paper.

This gives the following consequence.

\begin{ThmRef}{Main2}
If $\Gamma$ is a group with FDC and a finite classifying space $K(\Gamma,1)$ and if $R$ is a ring of finite global dimension, then the assembly maps $A_n \colon H_n (\Gamma, K(R)) \to K_n (R[\Gamma])$ are isomorphisms for all integers $n$.  In particular, the Whitehead group $\textit{Wh} (\Gamma)$ is trivial.
\end{ThmRef}


\begin{thebibliography}{99}

\bibitem{eAeM:11}
{E.~Aljadeff and E.~Meir},
\textit{Nilpotency of Bocksteins, Kropholler's hierarchy and a conjecture of Moore},  Adv. Math. \textbf{226} (2011), 4212--4224.

\bibitem{aB:03}
{A.C.~Bartels}, \textit{Squeezing and higher algebraic K-theory}
K-Theory, \textbf{28} (2003), 19--37.

\bibitem{pBmD:79}
P.H.~Berridge and M.J.~Dunwoody, \textit{Non-free projective modules for torsion-free groups}, J. London Math. Soc. \textbf{19} (1979), 193--215.

\bibitem{gCbG:03}
{G. Carlsson and B. Goldfarb}, \textit{On homological coherence of
discrete groups}, J.~Algebra {\bf 276} (2004), 502--514.

\bibitem{gCbG:04}
\bysame, \textit{The integral K-theoretic Novikov conjecture for
groups with finite asymptotic dimension}, Inventiones Math. {\bf 157} (2004), 405--418.

\bibitem{gCbG:11}
\bysame, \textit{Controlled algebraic $G$-theory, I},
J. Homotopy Relat. Struct. \textbf{6} (2011), 119--159.

\bibitem{gCbG:13}
\bysame, \textit{Algebraic $K$-theory of geometric groups}, preprint submitted for publication, 2013.
\texttt{arXiv:1305.3349}

\bibitem{aD:00}
{A.~Dranishnikov}, \textit{Asymptotic topology}, Russian Math. Surveys \textbf{55} (2000), 1085--1129.

\bibitem{aDmZ:12}
{A.~Dranishnikov and M.~Zarichnyi}, \textit{Asymptotic dimension, decomposition complexity, and Haver's property C}, preprint, 2012. \texttt{arXiv:1301.3484}

\bibitem{mD:72}
M.J.~Dunwoody, \textit{Relation modules}, Bull. London Math. Soc. \textbf{3} (1972), 151--155.

\bibitem{bG:03}
{B.~Grave},
\textit{A finitely presented group with infinite asymptotic dimension}, preprint, 2003.
available from \texttt{http:/$\!$/www.uni-math.gwdg.de/bgr/}

\bibitem{mG:93}
{M.~Gromov}, \textit{Asymptotic invariants of infinite groups}, in \textit{Geometric groups theory, Vol.2}, Cambridge U. Press (1993).

\bibitem{eGrTgY:12}
{E.~Guentner, R.~Tessera, and G.~Yu},
\textit{A notion of geometric complexity and its applications to topological rigidity}, Inventiones Math. \textbf{189} (2012), 315--357.

\bibitem{jHjJ:06}
J.~Harlander and J.A.~Jensen,
\textit{Exotic relation modules and homotopy types for certain 1-relator groups},
Alg. \& Geom. Top. \textbf{6} (2006), 2163--2173.

\bibitem{fJ:12}
F.E.A.~Johnson,
\textit{Syzygies and Homotopy Theory}, Springer-Verlag (2012).

\bibitem{pK:93}
{P.H.~Kropholler},
\textit{On groups of type $\mathrm{(FP)}_{\infty}$},
J. Pure Appl. Algebra \textbf{90} (1993), 55--67.

\bibitem{pK:99}
\bysame,
\textit{Modules possessing projective resolutions of finite type},
J. Algebra \textbf{216} (1999), 40--55.

\bibitem{jlL:76}
{J.-L.~Loday},
\textit{$K$-th\'eorie alg\'ebrique et repr\'esentations de groupes},
Ann. sci. \'Ecole Norm. Sup. \textbf{9} (1976), 309--377.

\bibitem{dRrTgY:11}
{D.~Ramras, R.~Tessera, and G.~Yu},
\textit{Finite decomposition complexity and the integral Novikov conjecture for higher algebraic K-theory},
to appear in J. Reine Angew. Math. (2013). \texttt{arXiv:1111.7022}

\bibitem{fW:78}
{F.~Waldhausen}, \textit{Algebraic $K$-theory of generalized free products}, Annals of Math. \textbf{108} (1978), 135--256.

\bibitem{gY:98}
{G.~Yu}, \textit{The Novikov conjecture for groups with finite asymptotic dimension}, Annals of Math. \textbf{147} (1998), 325--355.

\end{thebibliography}
\end{document}